\documentclass[12pt]{article}
\input{epsf.sty}
\usepackage{epsfig}
\usepackage{amsmath}
\usepackage{amssymb}
\usepackage{amsthm}
\usepackage{times}
\usepackage{indentfirst}
\newtheorem{mythm}{Theorem}

\newtheorem{mylem}{Lemma}
\newtheorem{myprop}{Proposition}
\newtheorem{mydef}{Definition}
\newtheorem{myrem}{Remark}
\def\sqr#1#2{{\vcenter{\vbox{\hrule height.#2pt \hbox{\vrule width.
#2pt height#1pt \kern#1pt \vrule} \hrule height.#2pt}}}}

\begin{document}
\begin{titlepage}
\begin{center}
{\large\bf Consensus Problems in Networks of Agents under Nonlinear
Protocols with Directed Interaction Topology}\footnote{This work was
supported by the National Science Foundation of China under Grant
No. 60574044, 60774074 and the Graduate Student Innovation
Foundation of Fudan University.}
\\[0.2in]
\begin{center}
Xiwei Liu\footnote{Email:xiwei.liu.fudan@gmail.com}, Tianping
Chen\footnote{These authors are with Lab. of Nonlinear Mathematics
Science, Institute of Mathematics, Fudan
University, Shanghai, 200433, P.R.China.\\
\indent ~~Corresponding author: Tianping Chen.
Email:tchen@fudan.edu.cn}
\end{center}

\end{center}
\begin{abstract}

The purpose of this short paper is to provide a theoretical analysis
for the consensus problem under nonlinear protocols. A main
contribution of this work is to generalize the previous consensus
problems under nonlinear protocols for networks with undirected
graphs to directed graphs (information flow). Our theoretical result
is that if the directed graph is strongly connected and the
nonlinear protocol is strictly increasing, then consensus can be
realized. Some simple examples are also provided to demonstrate the
validity of our theoretical result.

\end{abstract}
\end{titlepage}

\pagestyle{plain} \pagenumbering{arabic}

\section{Introduction}

In networks of dynamics agents, {\it ``consensus''} means that all
agents need to agree upon certain quantities of interest that
depend on their state. A {\it ``consensus protocol''} is an
interaction rule that specifies the information exchange between
an agent and all of its neighbors on the network, and enables the
network to achieve consensus via a process of distributed decision
making. Consensus problems (see \cite{SFM07}-\cite{RB05}) have a
long history in the field of computer science, particularly in
automata theory and distributed computation. Recently, distributed
coordination of networks of dynamic agents has attracted several
researchers from various disciplines of engineering and science
due to the broad applications of multi-agent systems in many
areas, such as collective behavior of flocks and swarms
\cite{GP03,RS06}, synchronization of coupled oscillators
\cite{Pik}-\cite{Chen2007}, and so on.

Until now, most papers in the literature mainly concern the
consensus problem under linear protocols, with the connection
topologies time-varying, state-dependent ( see
\cite{SFM07}-\cite{RB05}). Even in those papers investigating
nonlinear protocols, like \cite{SM03,SM04,Chen2007}, a strong
assumption on networks should be satisfied: the interaction topology
should be bidirectional. However, unidirectional communication is
important in practical applications and can be easily incorporated,
for example, via broadcasting. Also, sensed information flow which
plays a central role in schooling and flocking is typically not
bidirectional.

So, in this paper, we will look at the consensus problem in
networks of dynamic agents, described by ordinary differential
equations (ODE), under nonlinear protocols with directed topology.
This note can be regarded to extend consensus results under
undirected graphs in \cite{SM03,SM04} to the case of directed
graphs. Our approach is to model the communication topology as a
graph, then by merging spectral graph theory, matrix theory and
control theory, we can prove rigorously that if the directed graph
is strongly connected and the nonlinear protocol is strictly
increasing, then consensus problem can be realized.

An outline of this paper is as follows. In Section II, we define the
consensus problem on graphs. In Section III, we first define the
nonlinear protocol, then based on some lemmas of algebraic graph
theory and matrix theory, we obtain the main theoretical result. In
section IV, two simple examples are also provided to demonstrate the
effectiveness of our theoretical result. We conclude this paper in
Section V.

\section{Consensus problem on graphs}
\begin{mydef}(Weighted Directed Graph)
Let $\mathcal{G}=(\mathcal{V},\mathcal{E},\mathcal{A})$ be a
weighted digraph (or directed graph) with the set of nodes
$\mathcal{V}=\{v_1,\cdots,v_n\}$, set of edges
$\mathcal{E}\subseteq \mathcal{V}\times \mathcal{V}$, and a
weighted adjacency matrix $\mathcal{A}=(a_{ij})$ with nonnegative
adjacency elements $a_{ij}$. An edge of $\mathcal{G}$ is denoted
by $e_{ij}=(v_i,v_j)\in \mathcal{E}$, which means that node $v_i$
receives information from node $v_j$, and we assume that
$v_i\not=v_j$ for all $e_{ij}$, so the graph has no self-loops.
The adjacency elements associated with the edges of the graph are
positive, i.e., $e_{ij}\in \mathcal{E}\Longleftrightarrow
a_{ij}>0$. Moreover, we assume $a_{ii}=0$ for all $i\in
1,\cdots,n$. The set of neighbors of node $v_i$ is denoted by
$\mathcal{N}_i=\{v_j\in \mathcal{V}:(v_i,v_j\in\mathcal{E})\}$.
The corresponding {\it graph Laplacian} $L=(l_{ij})$ can be
defined as
\begin{eqnarray}
l_{ij}= \left\{\begin{array}{cc}
\sum_{k=1,k\not=i}^na_{ik},   &i=j\\
     -a_{ij},                 &i\not=j
\end{array}\right.
\label{Laplace}
\end{eqnarray}
\end{mydef}

\begin{mydef}(Strongly Connected Graph)
A path on a graph
$\mathcal{G}=(\mathcal{V},\mathcal{E},\mathcal{A})$ of length
$n^{\star} \leq n$ from $v_{i_0}$ to $v_{i_{n^{\star}}}$ is an
ordered set of distinct vertices
$\{v_{i_0},\cdots,v_{i_{n^{\star}}}\}$ such that
$(v_{i_{j-1}},v_{i_{{j}}})\in \mathcal{E}$, for all $j=1,\cdots,
n^{\star}$. A graph in which a path exists from every vertex to
every vertex is said to be {\it strongly connected (SC)}.
Obviously, irreducibility of the graph Laplacian for a graph can
imply its strong connectivity.
\end{mydef}

Without loss of generality, let $x_i\in R$ denote the value of
node $v_i$, $i=1,\cdots,n$. We refer to
$\mathcal{G}_x=(\mathcal{G},x)$ with $x=(x_1,\cdots,x_n)^T$ as a
network (or algebraic graph) with value $x$ and topology (or
information flow) $\mathcal{G}$. The value of a node might
represent physical quantities including attitude, position,
temperature, voltage, and so on.

\begin{mydef}(Consensus)
Consider a network of dynamic agents with $\dot{x}_i=u_i$ interested
in reaching a consensus via local communication with their neighbors
on a graph $\mathcal{G}_x$. By reaching a consensus, we mean
converging to a one-dimensional agreement space characterized by the
following equations:
\begin{eqnarray}
x_1=x_2=\cdots=x_n
\end{eqnarray}
This agreement space can be expressed as $x=\beta {\bf 1}$ where
${\bf 1}=(1,\cdots,1)^T$ and $\beta\in R$ is the collective
decision of the group of agents.
\end{mydef}

\begin{mylem}\label{dec}(See \cite{SFM07})
Suppose $L=(l_{ij})$ is a graph Laplacian of a bi-graph
$\mathcal{G}=(\mathcal{V},\mathcal{E},\mathcal{A})$ of $n$ nodes,
i.e., $l_{ij}=l_{ji}$, for any $i,j\in 1,\cdots,n$. The following
sum-of-squares (SOS) property holds, for any
$x=(x_1,\cdots,x_n)^T$,
\begin{eqnarray}
x^TLx=-\sum\limits_{j>i}l_{ij}(x_j-x_i)^2
\end{eqnarray}
\end{mylem}

\section{Consensus Analysis}
\subsection*{{\it A. Nonlinear consensus protocol}}
In this paper, we propose the following nonlinear consensus protocol
$h(\cdot): R\rightarrow R$ to solve consensus problems in a network
of continuous-time integrator agents with fix connection topology
$\mathcal{G}_x$:
\begin{eqnarray}
\dot{x}_i(t)=\sum\limits_{j\in
\mathcal{N}_i}{a}_{ij}\bigg(h(x_j(t))-h(x_i(t))\bigg),\qquad
i=1,\cdots,n \label{P1}
\end{eqnarray}
If $L=(l_{ij})$ is the corresponding graph Laplacian of
$\mathcal{G}_x$ defined in Definition 1, then the above equations
also can be rewritten as
\begin{eqnarray}
\dot{x}_i(t)=-\sum\limits_{j=1}^{n}l_{ij}h(x_j(t)),\qquad
i=1,\cdots,n \label{P1a}
\end{eqnarray}

Throughout this paper, we assume that $h(\cdot)$ is a strictly
increasing function. Without loss of generality, we assume $h(0)=0$.

\subsection*{{\it B. Algebraic graph theory and matrix theory}}
In this part, we introduce some basic concepts, notations and
lemmas in algebraic graph theory and matrix theory that will be
used throughout this paper.

\begin{mylem}
(Spectral localization. See \cite{SM04}) Let $\mathcal{G}$ be a
strongly connected digraph of $n$ nodes. Then $rank(L)=n-1$, and
all nontrivial eigenvalues of $L$ have positive real part.
\end{mylem}

\begin{myrem}
Lemma 2 holds under a weaker condition of existence of a directed
spanning tree for $\mathcal{G}$. $\mathcal{G}$ has a directed
spanning tree if there exists a node $r$ (root) such that all
other nodes can be linked to $r$ via a directed path (see relating
papers \cite{RB05,W05}). In fact, in digraphs with spanning tree
(leader-follower model), the root node is commonly known as a
leader, which does not receive any information from other nodes.
\end{myrem}

\begin{mylem}\label{left}
Assume $\mathcal{G}$ is a strongly connected digraph with graph
Laplacian $L$, then (\cite{H87})

1. ${\bf 1}=(1, 1, \cdots, 1)^T$ is the right eigenvector of $L$
corresponding to eigenvalue $0$ with multiplicity $1$;

2. Let ${\xi}=(\xi_{1}, \xi_2, \cdots, \xi_n)^T$ be the left
eigenvector of $L$ corresponding to the eigenvalue $0$. Then,
$\xi_{i}> 0$, $i= 1, 2, \cdots,n$; and its multiplicity is $1$. In
the following, we always assume $\sum_{i=1}^{n}\xi_{i}=1$.
\end{mylem}

\subsection*{{\it C. Main results}}
In this part, we will give a theorem, which shows that if the
directed graph is strongly connected and the nonlinear function is
strictly increasing, then the consensus problem can be realized.

\begin{mythm}
Suppose the digraph $\mathcal{G}_x$ is a strongly connected. $L$ is
the corresponding graph Laplacian in Definition 1. Then consensus
can be realized globally for all initial states by the nonlinear
protocol (\ref{P1a}) and the group decision is
$x_{\xi}=\sum_{i=1}^{n}\xi_ix_i(0)$, where
$\xi=(\xi_1,\cdots,\xi_n)^T$ is defined in Lemma 3.
\end{mythm}

Before the proof of Theorem 1, we introduce a reference node (or
virtue leader) $x_{\xi}(t)=\sum_{i=1}^n\xi_ix_i(t)$. It is clearly
\begin{eqnarray}
\dot{x}_{\xi}(t)=\sum\limits_{i=1}^n\xi_i\dot{x}_i(t)
=-\sum\limits_{i=1}^n\xi_i\sum\limits_{j=1}^n{l}_{ij}h(x_j(t))
=-\sum\limits_{j=1}^nh(x_j(t))\sum\limits_{i=1}^{n}\xi_i{a}_{ij}=0
\end{eqnarray}
Therefore, we obtain the following simple but useful proposition,
which plays an important role in the discussion of final group
decision.
\begin{myprop}
$x_{\xi}(t)$ is time-invariant for the network (\ref{P1a}), i.e.,
$x_{\xi}=\sum_{i=1}^n\xi_ix_i(0)=x_{\xi}(t)$ for all $t\geq 0$.
\end{myprop}

{\bf Proof of Theorem 1} Denote $x_{\xi}=\sum_{i=1}^n\xi_ix_i(0)$,
$x(t)=(x_1(t),\cdots,x_n(t))^T$, and $H(x(t))=(h(x_1(t)),
\cdots,h(x_n(t)))^T$. Then equations (\ref{P1a}) can be rewritten in
the compact form as
\begin{eqnarray}
\dot{x}(t)=-LH(x(t))\label{com}
\end{eqnarray}

Define a function as:
\begin{eqnarray}
V(x(t))=\sum\limits_{i=1}^n\xi_i\int_{0}^{x_i(t)}h(s)ds
\end{eqnarray}
Obviously, $V(x(t))\geq 0$ is radially unbounded, and $V(x(t))=0$
if and only if $x(t)=0$.

Denote $B=(b_{ij})=(\Xi L+L^T\Xi)/2$, where
$\Xi=\mathrm{diag}(\xi)$. It is easy to check that $B$ is a
symmetric matrix with zero row-sum, i.e., $B$ can be regarded as a
graph Laplacian of a bi-graph. Differentiating $V(x(t))$ and using
Lemma \ref{dec}, we have
\begin{eqnarray}
\dot{V}(x(t))&=&-\sum\limits_{i=1}^n\xi_ih(x_i(t))\sum\limits_{j\in
\mathcal{N}_i}{l}_{ij}h(x_j(t))\nonumber\\
&=&-H(x(t))^T{\Xi L}H(x(t))=-H(x(t))^TBH(x(t))\nonumber\\
&=&\sum\limits_{i>j}b_{ij}(h(x_i(t))-h(x_j(t)))^2\leq 0,~~~~
(\mathrm{since}~~b_{ij}\le 0) \label{dev}
\end{eqnarray}
Thus, $0\le V(x(t))\leq V(x(0))$, which implies $x(t)$ is bounded
for any $t\geq 0$. And the largest invariant subset ($\Omega$-limit
set) for the equations (\ref{P1a}) is
\begin{eqnarray}
\Omega=\{x:x_i(t)=x_j(t); i,j=1,\cdots,n\}\label{invariance}
\end{eqnarray}

Now, we claim that for all $i=1,\cdots,n$,
$\lim\limits_{t\rightarrow\infty}x_{i}(t)=x_{\xi}$.

In fact, if $t_{m}\rightarrow\infty$ and for all $i,j=1,\cdots,n$,
$x_{i}(t_{m})\rightarrow\beta$, then,
\begin{eqnarray}
x_{\xi}=
\lim\limits_{t\rightarrow\infty}\sum\limits_{i=1}^{n}\xi_ix_i(t_{m})
=\sum\limits_{i=1}^{n}\xi_i\beta=\beta
\end{eqnarray}
which means that $x_{\xi}$ is the group decision of the consensus
problem. Theorem 1 is proved completely.

\begin{myrem}
Let $h(x_i(t))=\alpha x_i(t)$ with $\alpha>0$, The nonlinear
function $h(\cdot)$ becomes a linear function. therefore, Theorem 1
can be regarded as a generalization of the consensus problem under
linear protocols. It also give a simple proof for the consensus
problem under linear protocols, too.
\end{myrem}

\begin{myrem}
Assume $(h(w_1)-h(w_2))/(w_1-w_2)\geq \alpha$ holds for $\alpha>0$
and any $w_1\not=w_2\in R$. In this case, we have
\begin{eqnarray}
\dot{V}(x(t))=\sum\limits_{i>j}b_{ij}(h(x_i(t))-h(x_j(t)))^2
\leq\alpha^2\sum\limits_{i>j}b_{ij}(x_i(t)-x_j(t))^2
\end{eqnarray}
and the consensus problem will be realized exponentially. Moreover,
it seems that the nonlinear protocol can be realized faster than
that under the linear protocol $h(w)=\alpha w$. Therefore, nonlinear
protocols can be applied to calculate the average value of
large-scale networks more effectively.
\end{myrem}

\section{Numerical examples}
In this section, we give two numerical simulations to verify the
validity of our theory.

Consider a network of a strongly connected digraph $\mathcal{G}_x$
with $3$ agents
\begin{eqnarray*}
\dot{x}_i(t)=-\sum\limits_{j=1}^3{l}_{ij}h(x_j(t)),\qquad i=1,2,3
\label{cex}
\end{eqnarray*}
where $x_i(t)\in R$ and  the Laplacian of $\mathcal{G}_x$ is
\begin{eqnarray}
L=\left(\begin{array}{ccc}
2  &  -1&  -1\\
0  &   1&  -1\\
-1 &   0&   1
\end{array}\right)\label{matrix}
\end{eqnarray}
Its left eigenvector with eigenvalue $0$ is $\xi=(1/4,1/4,1/2)^T$.

{\bf Example 1:}\quad In this simulation, the nonlinear protocol is
assumed as $h(x_i(t))=\alpha x_i(t)+\sin(x_i(t)),~ i=1,2,3$. The
initial value is taken as: $(x_1(0),x_2(0),x_3(0))=(1,2,3)$.

\begin{figure}
\centering
\includegraphics[width=0.8\textwidth, height=0.8\textheight]{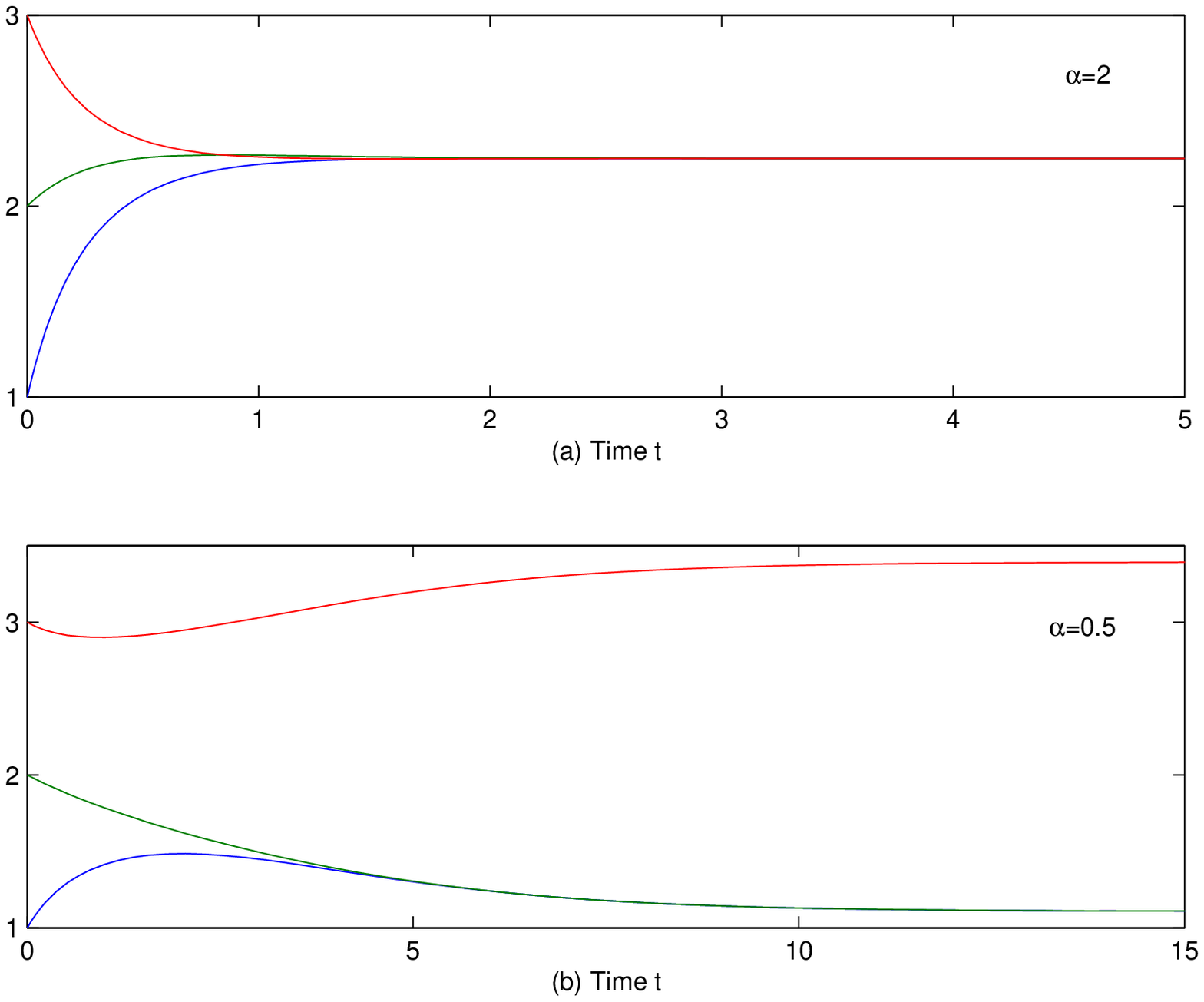}
\caption{Consensus problem of (\ref{cex}) under different
nonlinear functions} \label{cf}
\end{figure}

Case 1. $\alpha=2$. In this case, then $h^{\prime}(\cdot)\geq 1$. By
Theorem 1, consensus of (\ref{cex}) can be realized, and the
decision value is $\sum_{i=1}^3\xi_ix_i(0)={1}/{4}+2/4+3/2=2.25$,
see Figure \ref{cf}(a);

Case 2. $\alpha=0.5$. In this case, $h(\cdot)$ is not an increasing
function, and consensus of (\ref{cex}) may not be realized, see
Figure \ref{cf}(b).

{\bf Example 2:}\quad In this simulation, we choose two protocols.
One is the nonlinear protocol
\begin{eqnarray}
h^{\star}(x_i)=\left\{\begin{array}{ccl}
x_i^2           & ;&\mathrm{if}\quad x_i>1\\
\sqrt{x_i}      & ;&\mathrm{if}\quad 0<x_i\leq 1\\
-\sqrt{-x_i}    & ;&\mathrm{if}\quad -1<x_i\leq 0\\
-x_i^2          & ;&\mathrm{if}\quad x_i\leq-1
\end{array}\right.\label{non}
\end{eqnarray}
The other is the linear protocol
\begin{eqnarray}
h(x_i)=x_i/2\qquad i=1,2,3 \label{lin}
\end{eqnarray}

Simple calculations show that the derivative of $h^{\star}(\cdot)$
is no less than $1/2$. The consensus problem under the nonlinear
protocol (\ref{non}) can be realized faster than that under the
linear one (\ref{lin}).

\begin{figure}
\centering
\includegraphics[width=0.8\textwidth, height=0.4\textheight]{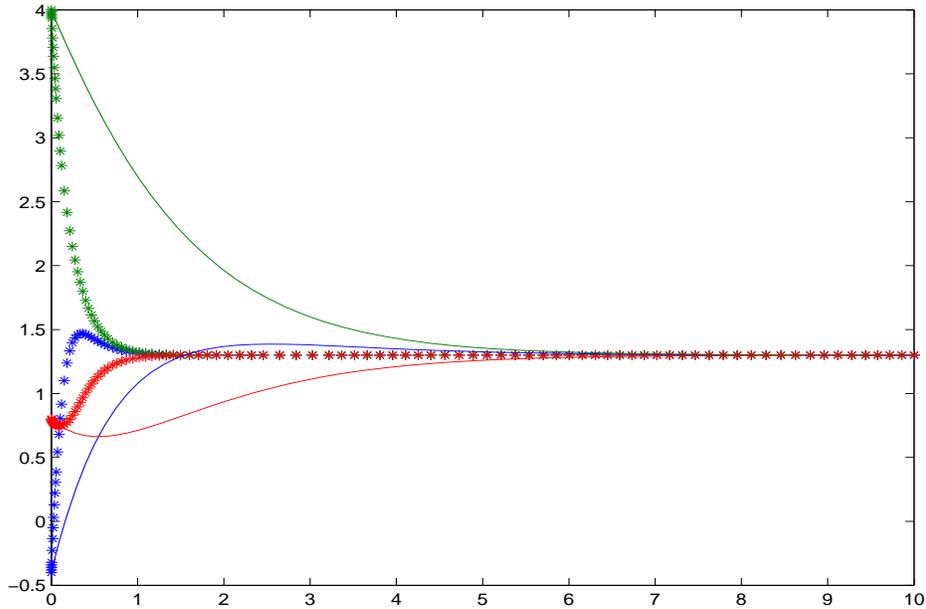}
\caption{Consensus problem of (\ref{cex}) under nonlinear and linear
protocols} \label{cf2}
\end{figure}

In Figure \ref{cf2}, the dynamical behavior of the network
(\ref{cex}) under the nonlinear protocol $h^{*}(x)$ defined in
(\ref{non}) is displayed by line with star. Instead, for the linear
protocol $h(x)$ defined in (\ref{lin}), it is displayed by line
without star. The initial value is chosen as:
$(x_1(0),x_2(0),x_3(0))=(-0.4,4,0.8)$. Simulations do show that
consensus under the nonlinear protocol (\ref{non}) is much faster
than that under linear protocol (\ref{lin}).

\section{Conclusions}
In this paper, we investigate the consensus problem under nonlinear
protocols. We generalize the results for undirected graphs to
directed graphs. Moreover, our model can also be regarded as the
generalization of consensus problem under linear protocols to
nonlinear protocols. All the existing results with respect to
consensus under linear protocols with directed/undirected graph and
consensus under nonlinear protocols with undirected graph can be
easily obtained by our approach. The convergence analysis is
presented rigorously, based on tools from algebraic graph theory,
matrix theory, and control theory. Two simple examples are provided
to show the effectiveness of the theoretical result.

\end{document}